\documentclass[leqno,a4paper,11pt]{article}

\usepackage{hyperref}
\usepackage{amsthm}
\usepackage{amsmath}
\usepackage{amssymb}
\usepackage{graphicx}
\usepackage{centernot}
\usepackage[xcolor]{changebar}
\usepackage{enumitem}

\newcommand{\Con}{\ensuremath{\mathcal{C}}}

\newcommand{\mb}[1]{\ensuremath{\mathbb{#1}}}
\newcommand{\N}{\mb{N}}

\newcommand{\R}{\mb{R}}

\newfont{\bl}{msbm10 scaled \magstep2}

\newcommand{\beq}{\begin{equation}}
\newcommand{\eeq}{\end{equation}}

\newcommand{\notmid}{\mid\kern-0.5em\not\kern0.5em}

\theoremstyle{theorem}
\newtheorem{thm}{Theorem}[section]
\newtheorem{lem}[thm]{Lemma}
\newtheorem{prop}[thm]{Proposition}
\newtheorem{cor}[thm]{Corollary}
\newtheorem{defi}[thm]{Definition}

\theoremstyle{definition}
\newtheorem{rem}[thm]{Remark}
\newenvironment{pr}{\begin{proof}[\textbf{Proof:}] \ }{\end{proof}}
\newtheorem{ex}[thm]{Example}

\let\oldmarginpar\marginpar
\renewcommand\marginpar[1]{\-\oldmarginpar[\raggedleft\footnotesize #1]%
{\raggedright\footnotesize #1}}

\newcommand{\LLS}{Lorentzian length space }
\newcommand{\LLSn}{Lorentzian length space}
\newcommand{\LpLS}{Lorentzian pre-length space }
\newcommand{\LpLSn}{Lorentzian pre-length space}
\newcommand{\Xll}{$(X,d,\ll,\leq,\tau)$ }
\newcommand{\Xlln}{$(X,d,\ll,\leq,\tau)$}
\newcommand{\tXll}{$(\tilde X,\tilde d,\tilde \ll,\tilde \leq,\tilde \tau)$ }

\newcommand{\Tau}{\mathcal{T}}

\title{Inextendibility of spacetimes and Lorentzian length spaces} \date{}
 \author{James D.E.\ Grant\thanks{{\tt j.grant@surrey.ac.uk}, Department of Mathematics,
University of Surrey, UK.}, Michael Kunzinger\thanks{{\tt michael.kunzinger@univie.ac.at}, Faculty of 
Mathematics, University of Vienna, Austria.},
Clemens S\"amann\thanks{{\tt clemens.saemann@univie.ac.at}, Faculty of Mathematics, University of Vienna,
 Austria.}}

\begin{document}
 \maketitle

 \begin{abstract}
We study the low-regularity (in-)extendibility of spacetimes within the synthetic-geometric framework of Lorentzian length spaces 
developed in~\cite{KS:18}. To this end, we introduce appropriate notions of geodesics and timelike geodesic completeness and 
prove a general inextendibility result. Our results shed new light on recent analytic work in this direction and, for the 
first time, relate low-regularity inextendibility to (synthetic) curvature blow-up. 
\vskip 1em

\noindent
\emph{Keywords:} Length spaces, Lorentzian length spaces, causality theory, synthetic curvature bounds, triangle comparison,
metric geometry, inextendibility
\medskip

\noindent
\emph{MSC2010:} 53C23, 
53C50, 
53B30, 
53C80, 
83C75  
 \end{abstract}

\section{Introduction}\label{sec-intro}

One can distinguish between two main lines of research in low-regularity geometry. One approach is analytical, where one 
lowers the differentiability assumptions on, for example, (pseudo-)Riemannian metrics below the level where curvature can be 
classically defined. For example, one can study geometrical properties of (pseudo-)Riemannian metrics that have regularity  
$C^0$, $C^{0, \alpha}$ or $C^{1, 1}$, etc., or so-called ``Geroch--Traschen'' metrics, for which the Christoffel symbols are 
$L^2_{\mathrm{loc}}$, and the curvature is well-defined as a distribution~\cite{GT, LM:07, SV:09}. The other approach to 
studying low-regularity geometries is by ``synthetic'' or metric space methods. Here, curvature bounds for Alexandrov spaces 
and CAT$(k)$ spaces are defined in terms of comparison properties of geodesic triangles. 

In the context of low-regularity Riemannian geometry, examples of a result of an analytical nature would be DeTurck and 
Kazdan's study concerning harmonic coordinates~\cite{DK}, Taylor's results on regularity of isometries~\cite{T:06} and 
Lytchak and Yaman's result~\cite{LY} that minimising curves for $C^{0, \alpha}$ Riemannian manifolds are $C^{1, \beta}$ 
curves, where $\beta = \frac{\alpha}{2-\alpha}$. Examples in this direction in the Lorentzian setting are the positive mass 
theorem for distributional curvature \cite{LL:15, GT:14}, work on cone structures \cite{FS:12, Min:18, BS:18} and the recent 
work of extending the classical singularity theorems to $\Con^{1,1}$-regularity \cite{KSSV:15, KSV:15, GGKS:18}, which in 
turn builds on previous results in low regularity Lorentzian geometry and causality \cite{CG:12, Min:15, KSS:14, KSSV:14, 
Sae:16}.

In the synthetic direction, the theory of Alexandrov spaces with curvature bounded above and/or below is well-developed as an 
appropriate generalisation of Riemannian geometry with sectional curvature bounds (see, for instance, \cite{BBI:01, BH:99, 
Pap:14}), and the work of Lott--Villani--Sturm gives a generalisation of the notion of a Riemannian metric with lower bound 
on the Ricci curvature to metric measure spaces~\cite{LV, S1, S2}.

In this paper, we will concentrate on a generalisation of Lorentzian geometry suitable for the low-regularity setting. More 
precisely, we shall be interested in the problem of finding low-regularity extensions of spacetimes. Concerning this 
question, approached from the analytical side, several fundamental contributions have appeared recently. Of particular 
relevance to us, Sbierski has shown the $C^0$-inextendibility of the Schwarzschild solution~\cite{Sbi:18}. Building upon 
Sbierski's work, Galloway, Ling and Sbierski established that global hyperbolicity combined with timelike geodesic 
completeness implies $C^0$-inextendibility. Further developments in this direction are due to Galloway--Ling and Graf--Ling 
(see below). In a related direction, Dafermos and Luk have recently shown
$C^0$-extendibility of the interior of the 
Kerr solution~\cite{DL:17}.

In this paper, we will concentrate on the synthetic-geometrical approach to extendibility. In~\cite{KS:18}, the theory of 
Lorentzian length spaces has been developed, which will form the framework of the present work. In this more axiomatic 
approach, there is a notion of a geodesic (as a locally length maximising curve), which is not available in the 
more analytical direction of research. Therefore, it is possible to mimic the classical proof that geodesic completeness implies 
inextendibility (see, for example,~\cite[Prop.~6.16]{BEE:96}). Moreover, within this picture, it becomes clear precisely what 
minimal geometric properties are underlying certain analytical extension results. In particular, for the first time, our 
approach allows us to directly relate low-regularity inextendibility with (synthetic) curvature blow-up. Such a result does 
not appear to be feasible in a purely analytical approach, due to the lack of a notion of a curvature for the extended 
spacetime. 

An additional advantage of our synthetic approach is that there is no requirement for the introduction of coordinate 
systems, and regularity conditions (such as existence of smooth structures, or a certain level of differentiability) never 
arise. In this regard, it should perhaps be noted that in the analytical work on low-regularity extensions, one has to carry 
out standard geometrical constructions on the original manifold. As such, even though one works in a coordinate chart of the 
extended manifold in which the metric is merely continuous, the metric on the intersection of the original manifold with the 
coordinate chart must be $C^2$-regular.%
\footnote{I.e.\ one implicitly must assume that the metric is smooth on $\iota(M)$ in the coordinate chart on $\tilde{M}$ in 
which the metric is just continuous.} 
One could compare this situation with, for example, the fact that the Nash--Kuiper theorem~\cite{Nash, K:56} 
implies that the flat metric on $T^2$ can be induced from a $C^1$ map $T^2 \to \R^3$.%
\footnote{See, for instance,~\cite{BJLT} for an illustration of this example.} 
In the coordinate system in which the map is $C^1$, the induced metric will be merely $C^0$, even though we know that there exists a coordinate system in which the metric is smooth.  As such, one could consider a more general notion of $C^0$ extensions of spacetimes, where one allows the regularity of the metric on the original manifold to drop. On the contrary, in our approach, such issues never arise. In fact, the extensions that we consider need not even be manifolds. 

Our main references for Lorentzian geometry and causality theory are \cite{ONe:83,BEE:96,MS:08,Chr:11},
as well as \cite{CG:12} for the case of continuous Lorentzian metrics. 

The plan of the paper is as follows: In Section \ref{sec:lls} we briefly recall some main concepts
and results on Lorentzian length spaces. Section \ref{sec:extensions} introduces extensions of
Lorentzian (pre-)length spaces, relates them to extensions of spacetimes and shows that the 
future or past boundary of an extension is non-empty. In Section \ref{sec:geodesics} we
define geodesics in the synthetic setting and show that this notion reduces precisely to
that of pregeodesics for spacetimes. We also demonstrate that, as in the smooth case, extendibility as a geodesic
is equivalent to continuous extendibility. In Section \ref{sec:TC} we define an analogue of
timelike completeness: a Lorentzian pre-length space is said to have property (TC) if all 
inextendible timelike geodesics have infinite length. This is the key property on which our main
inextendibility result (Theorem \ref{thm-inext-lls-lls}) rests. We then establish connections between inextendibility
and the occurrence of synthetic causal curvature singularities. Finally, in Section \ref{sec:relation}
we relate the results of the present work to the recent advances in the study of 
the low regularity inextendibility
of spacetimes.

\section{A short introduction to \LLSn s}\label{sec:lls}
Here we briefly recall some basic notions and results
from the theory of Lorentzian length spaces, following \cite{KS:18}, to which we
refer for further details and proofs.

A set $X$ endowed with a preorder $\leq$ and a transitive relation  $\ll$ 
contained in $\leq$ is called a \emph{causal space}. We write $x<y$ if $x\leq y$ and $x\neq y$.
If $x\ll y$ respectively $x\le y$ we call $x$ and $y$ timelike respectively causally related. 
Chronological and causal futures and pasts $I^\pm(x)$, $J^\pm(x)$ of a point $x$
are then defined in the usual manner based on these relations.

If $X$ is, in addition, equipped with a metric  $d$
and a lower semicontinuous map
$\tau \colon X\times X \to [0, \infty]$ that satisfies the reverse triangle inequality
$\tau(x,z)\geq \tau(x,y) + \tau(y,z)$ (for all $x\leq y\leq z$), as well as $\tau(x,y)=0$ if $x\nleq y$ and $\tau(x,y)>0 
\Leftrightarrow x\ll y$, then \Xll is called a \emph{Lorentzian pre-length space\/} and $\tau$ is called 
the \emph{time separation function\/} of $X$. Note that lower semicontinuity of $\tau$ implies that $I^{\pm}(x)$ is open, for any $x \in X$. 

A non-constant curve $\gamma \colon I\rightarrow X$ ($I$ an interval) is called 
(future-directed) \emph{causal (timelike)} if $\gamma$ is locally Lipschitz continuous and if for all 
$t_1,t_2\in I$ with $t_1<t_2$ we have $\gamma(t_1)\leq\gamma(t_2)$ ($\gamma(t_1)\ll\gamma(t_2)$). It
is called \emph{null\/} if, in addition to being causal, no two points on the curve are related with respect 
to $\ll$. For strongly causal continuous Lorentzian metrics, this notion of causality
coincides with the usual one (\cite[Prop.\ 5.9]{KS:18}). In analogy to the theory
of metric length spaces, the length of a causal curve is defined via the time separation
function:  For $\gamma \colon [a,b]\rightarrow X$ future-directed causal we set
\[
L_\tau(\gamma):=
\inf\Big\{\sum_{i=0}^{N-1} \tau(\gamma(t_i),\gamma(t_{i+1})): a=t_0<t_1<\ldots<t_N=b,\ N\in\N\Big\}.
\ 
\]
If the interval is (half-)open, say $I=[a,b)$, then the infimum is taken over all partitions
with $a=t_0<t_1<\ldots<t_N<b$, and similarly for the other cases. 
For smooth and strongly causal spacetimes $(M,g)$ this notion of length coincides with the usual one:
$L_\tau(\gamma)=L_g(\gamma)$ (\cite[Prop.\ 2.32]{KS:18}).
A future-directed causal curve $\gamma \colon [a,b]\rightarrow X$ is 
\emph{maximal\/} if it realizes the time separation, i.e., if $L_\tau(\gamma) = \tau(\gamma(a),\gamma(b))$.

Standard causality conditions (chronology, (strong) causality, global hyperbolicity, \dots) can also be imposed
on Lorentzian pre-length spaces, and substantial parts of the causal ladder (\cite{MS:08}) continue
to hold in this general setting. A \LpLS $X$ is called \emph{causally path connected\/} if for 
all $x,y\in X$ with $x\ll y$ (respectively $x<y$) there is a future-directed 
timelike (respectively causal) curve from $x$ to $y$. A neighbourhood $U$ of $x$ is called 
\emph{causally closed\/} if the relation
$\leq$ is closed in $\bar{U}\times\bar{U}$, and $X$ itself is called 
\emph{locally causally closed\/} if every point has a causally closed neighbourhood. 

A key technical tool in smooth semi-Riemannian geometry is the existence of convex neighbourhoods, 
in which the causality is particularly simple and where one has a complete description of length-maximising
curves. The analogue of this notion in the present context is the following: 
A \LpLS $X$ is called \emph{localisable\/} if any $x\in X$ has an open, so-called \emph{localising\/} %
neighbourhood $\Omega_x$ such that:
\begin{enumerate}[label=(\roman*)]
	\item \label{def-loc-LpLS-cau-comp} The $d$-length of all causal curves contained in $\Omega_x$
	is uniformly bounded.
	\item \label{def-loc-LpLS-om-con} There is a continuous map $\omega_x\colon \Omega_x \times \Omega_x\rightarrow 
	[0,\infty)$ such that 
	$(\Omega_x, d\rvert_{\Omega_x\times\Omega_x},$ $\ll\rvert_{\Omega_x\times \Omega_x},\leq\rvert_{\Omega_x\times\Omega_x}, 
	\omega_x)$ is a Lorentzian pre-length space,  and for every $y\in\Omega_x$ we have 
	$I^\pm(y)\cap\Omega_x\neq\emptyset$.
	\item \label{def-loc-LpLS-max-cc} For all $p,q\in \Omega_x$ with $p<q$ there is a future-directed causal curve 
	$\gamma_{p,q}$ from $p$ to $q$ that is 
	maximal in $\Omega_x$ and satisfies	$L_\tau(\gamma_{p,q}) = \omega_x(p,q) \leq \tau(p,q)$.
\end{enumerate}
If, in addition, the neighbourhoods $\Omega_x$ can be chosen such that
\begin{enumerate}
	\item[(iv)]\label{def-loc-LpLS-4} Whenever $p,q\in\Omega_x$ satisfy $p\ll q$ then $\gamma_{p,q}$ is timelike 
	and strictly longer than any future-directed 
	causal curve in $\Omega_x$ from $p$ to $q$ that contains a null segment,
\end{enumerate}
then \Xll is called {\em regularly localisable}.

Lorentzian length spaces are close analogues of metric length spaces in the sense that the time separation function
can be calculated from the length of causal curves connecting causally related points. Precisely, a 
locally causally closed, causally path connected and localisable \LpLS is called a \emph{\LLSn\/} if
$\tau = \mathcal{T}$, where for any  $x,y\in X$ we set
\begin{equation*}
\mathcal{T}(x,y):= \sup\{L_\tau(\gamma):\gamma \text{ future-directed causal from }x \text{ to } y\}\,, 
\end{equation*}
if the set of future-directed causal curves from $x$ to $y$ is not empty. Otherwise let $\mathcal{T}(x,y):=0$. 
If, in addition, $X$ is regularly
localisable, then it is called a regular Lorentzian length space.

Any smooth strongly causal spacetime is an example of a regular Lorentz\-ian length space (with metric $d = d^h$ induced by any Riemannian metric $h$ on the spacetime). 
More generally,
any spacetime with a continuous, strongly causal and causally plain metric (see the remark preceding Corollary~\ref{cor-geo-compl-inext-lls} below) is a (strongly)
localisable Lorentzian length space. Further examples are provided by certain Lorentz-Finsler spaces
in the sense of \cite{Min:18} or, for the non-manifold setting, causal Fermion systems \cite{Fin:16,Fin:17}.

The final concept from the theory of Lorentzian length spaces we are going to require below is
that of synthetic curvature bounds, based on triangle comparison. We will confine ourselves to
causal triangle comparison here, as this is the only one we are going to employ. By an
\emph{admissible causal geodesic triangle\/} we mean a triple $(x,y,z)\in X^3$ with 
$x\ll y \le z$ or $x \le y \ll z$ such that $\tau(x,z)<\infty$ and such that the sides (if non-trivial) 
are realized by future-directed causal curves. Curvature bounds are formulated by comparing
such triangles with triangles of the same side lengths in one of the Lorentzian model spaces
$M_K$ of constant sectional curvature. Here, 
\begin{equation}\label{eq:model_spaces}
M_K = \left\{ \begin{array}{ll}
\tilde S^2_1(r) & K=\frac{1}{r^2}\\
\R^2_1 & K=0\\
\tilde H^2_1(r) & K= -\frac{1}{r^2}.
\end{array}
\right.
\end{equation}
where $\tilde S^2_1(r)$ is the simply connected covering manifold of the two-dimensional Lorentzian pseudosphere
$S^2_1(r)$ (i.e., de-Sitter space), $\R^2_1$ is two-dimensional Minkowski space, and $\tilde H^2_1(r)$ is the simply connected covering manifold
of the two-dimensional Lorentzian pseudohyperbolic space (i.e., anti-de-Sitter space) . In order to guarantee the existence of comparison 
triangles in one of the model spaces, one needs to impose size restrictions of the following kind:
Given $K\in \R$, let $(a,b,c)\in \R_+^3$ with $c\ge a+b$. 
If $c=a+b$, then let $c<\frac{\pi}{\sqrt{K}}$. (Here, $\frac{\pi}{\sqrt{K}}:=\infty$ 
if $K\le 0$). 
Otherwise, if $K<0$ then assume $c<\frac{\pi}{\sqrt{-K}}$.  Then $(a,b,c)$ is  said to {\em satisfy timelike size bounds} for $K$. These bounds ensure the existence of comparison triangles in the corresponding model space.

Using this terminology, a Lorentzian pre-length space \Xll is said to have causal curvature bounded below (above) by $K\in\R$ 
if every point in $X$ has a neighbourhood $U$ such that:
\begin{enumerate}[label=(\roman*)] 
	\item $\tau|_{U\times U}$ is finite and continuous.
	\item Whenever $x,y \in U$ with $x < y$, there exists a causal curve $\alpha$ in $U$ with $L_\tau(\alpha) = \tau(x,y)$.
	\item If $(x,y,z)$ is an admissible causal geodesic triangle in $U$, realized by maximal causal curves 
	(or a constant curve, respectively) $\alpha, \beta, \gamma$
	whose side lengths satisfy timelike size bounds for $K$, and if $(\bar{x},\bar{y},\bar{z})$ is a comparison triangle of 
	$(x,y,z)$ in $M_K$ realized by causal geodesics (or a constant curve) $\bar\alpha$, $\bar{\beta}$, $\bar{\gamma}$,
	then whenever $p$, $q$ are points on the timelike sides of $(x,y,z)$ and $\bar p$, $\bar q$ are corresponding
	points of the timelike sides of $(\bar{x},\bar{y},\bar{z})$,
	we have $\tau(p,q)\le \bar{\tau}(\bar p, \bar q)$ $($respectively $\tau(p,q)\ge \bar{\tau}(\bar p, \bar q))$.
\end{enumerate}
Such a neighbourhood $U$ is called a \emph{comparison neighbourhood with respect to $M_K$}.

\section{Extensions}\label{sec:extensions}
We start the main part of our work by defining the notion of an \emph{extension\/} of a \LpLSn, requiring only conditions that are natural within our setting.%
This concept is fully compatible with the usual notion of extension for spacetimes, see Proposition~\ref{prop-ext-ext}.

\begin{defi}\label{def-ext}
 Let \Xll be a \LpLSn . A \LpLS \tXll is called an \emph{extension} of \Xll if
 \begin{enumerate}[label=(\roman*)]
  \item \label{def-ext-conn} the metric space $(\tilde X, \tilde d)$ is connected,
  \item there exists an isometry $\iota\colon (X,d)\rightarrow (\tilde X, \tilde d)$ of metric spaces,
  \item the image $\iota(X)$ is a proper, open subset of $\tilde X$,
  \item \label{def-ext-cr} $\iota$ preserves $\ll$ and $\leq$, i.e., $\forall x,y\in X$: if $x\leq y$ then $\iota(x)\ 
\tilde{\leq}\ \iota(y)$ and if $x\ll y$ then $\iota(x)\ \tilde\ll\ \iota(y)$, and 
  \item \label{def-ext-tau} a curve $\gamma \colon I \to X$ is timelike (respectively causal) if and only if
  $\iota\circ \gamma$ is timelike (respectively causal) in $(\tilde X,\tilde d,\tilde \ll,\tilde \le,\tilde \tau)$.
  Furthermore, $\iota$ preserves $\tau$-lengths, i.e., for any $\leq$-causal curve $\gamma\colon I\rightarrow 
X$ we have
  \begin{equation}\label{eq-def-ext-tau}
   L_\tau(\gamma) = L_{\tilde\tau}(\iota\circ\gamma)\,.
  \end{equation}
 \end{enumerate}
 In this case \Xll is called \emph{extendible}. If no extension exists, then \Xll is called \emph{inextendible} (as a 
\LpLSn).
\end{defi}

\begin{rem}
Of course, this definition also applies to \LLSn s, i.e., a \LLS is \emph{extendible} if there is a \LLS \tXll and  
$\iota\colon (X,d)\rightarrow (\tilde X, \tilde d)$ with the above properties \ref{def-ext-conn}-\ref{def-ext-tau}. In 
this case conditions \ref{def-ext-cr} and \ref{def-ext-tau} slightly simplify. 
\end{rem}

\begin{lem}\label{lem-tau-ttau}
 Let \tXll be an extension of \Xlln, where both are \LLSn s. Then $\tilde\tau \circ (\iota\times\iota) \geq \tau$.
\end{lem}
\begin{pr}
 Let $p,q\in X$ with $\tau(p,q)>0$ (if $\tau(p,q)=0$ there is nothing to do). Let $\gamma$ be a future 
directed $\leq$-causal curve from $p$ to $q$ (which exists due to $p\leq q$ and the causal path-connectedness of $X$).
Then $\iota\circ\gamma$ is $\tilde\leq$-causal and $L_\tau(\gamma)=L_{\tilde\tau}(\iota\circ\gamma)\leq \tilde\Tau(\iota(p),\iota(q)) 
= \tilde\tau(\iota(p),\iota(q))$. Taking the supremum over all future-directed $\leq$-causal curves from $p$ to $q$ we get $\Tau(p,q)\leq 
\tilde\tau(\iota(p),\iota(q))$ and since $\Tau=\tau$ the claim follows.
\end{pr}

The following lemma shows that condition \ref{def-ext-tau} of Definition \ref{def-ext} required of an extension is in 
fact not too strong. Moreover, it demonstrates that for smooth strongly causal spacetimes the time separation function 
determines the metric completely.
\begin{lem}
 Let $(M,g)$ and $(\tilde M,\tilde g)$ be smooth spacetimes (of the same dimension) with time separation functions $\tau$ 
and $\tilde\tau$, respectively. Let $(M,g)$ be strongly causal and let $\iota\colon M\rightarrow \tilde M$ be onto. Then 
$\iota$ is an isometry if and only if $\iota$ preserves 
causal curves and their lengths, i.e., a curve $\gamma$ is causal in $M$
if and only if $\iota\circ\gamma$ is causal in $\tilde M$ and for such
curves, $L_g(\gamma) = 
L_{\tilde{g}}(\iota\circ\gamma)$.
\end{lem}
\begin{pr}
 It is a classical result that goes back to Hawking, King and McCarthy \cite{HKM:76} (cf.\ \cite[Prop.\ 
3.34]{MS:08} or \cite[Thm.\ 4.17]{BEE:96}) that $\iota$ is an isometry if and only if it preserves $\tau$.
By definition of the time separation functions in spacetimes,
this latter condition is, in turn, implied by $\iota$ preserving the $g$-lengths 
of causal curves.
\end{pr}

Furthermore, in the case of spacetimes the above result implies that there is no difference between an 
extension in our sense, and in the usual sense of an isometric embedding (cf.\ \cite[Def.\ 2.15]{Sbi:18}.%
To be precise, we have the following result:

\begin{prop}\label{prop-ext-ext}
  Let $(M,g)$ and $(\tilde M, \tilde g)$ be smooth, strongly causal spacetimes (of the same dimension) and let $\iota\colon 
M\rightarrow \tilde M$ be a map such that $\iota(M)\subset \tilde M$. Then the induced \LLS of $(\tilde M, 
\tilde g)$ extends the one coming from $(M,g)$ via $\iota$ if and only if $\iota$ is a (smooth) isometric embedding.
\end{prop}
\begin{pr}
We start with the following observation: Let $\tilde h$ be any Riemannian metric on $\tilde M$ with induced metric 
$d^{\tilde h}$. This fixes the induced \LLS in the following sense: Any other Riemannian metric on $\tilde M$ also 
induces the manifold topology and the notion of locally Lipschitz continuous curves is preserved (cf.\ \cite[Prop.\ 
2.3.1]{Chr:11}), thus fixing the spacetime $(\tilde M, \tilde g)$ and any Riemannian background metric determines 
the resulting \LLSn.

Assume that $(\tilde M,d^{\tilde h},\tilde\ll,\tilde\leq,\tilde \tau)$ extends $(M,d^h,\ll,\leq,\tau)$ via $\iota$. As 
$\iota(M)$ is an open and connected subset of $\tilde M$ we consider the spacetime $(\hat M,\hat g):= (\iota(M),\tilde 
g\rvert_{\iota(M)})$ with its time separation function $\hat \tau$. This means that
\begin{equation*}
 \hat\tau(\tilde p,\tilde q) = 
\sup\{L_{\tilde g}(\tilde \gamma): \tilde \gamma \text{ f.d.\ causal curve from } \tilde p\text{ to }\tilde q \text{ with } 
\textrm{image}(\tilde\gamma)\subseteq\iota(M)\}\,.
\end{equation*}
By Definition \ref{def-ext},\ref{def-ext-tau} a curve $\gamma\colon I\rightarrow M$ is causal if and only if 
$\iota\circ\gamma\colon I 
\rightarrow \hat M$ is causal in $(\hat M,\hat g)$. This together with \eqref{eq-def-ext-tau} and \cite[Prop.\ 2.32]{KS:18}
implies that $\iota$ preserves $\hat\tau$, i.e.,
\begin{equation*}
 \tau(p,q)=\hat\tau(\iota(p),\iota(q))\quad \forall p,q\in M\,.
\end{equation*}
Thus by \cite[Prop.\ 3.34]{MS:08} $\iota$ is an isometry $(M,g)\rightarrow (\hat M,\hat g)$.

\medskip

For the converse assume that $\iota$ is a smooth isometric embedding. Then we check points 
\ref{def-ext-conn}-\ref{def-ext-tau} of Definition \ref{def-ext}. As $\tilde M$ is connected by 
assumption,
the first point follows. Pulling back $\tilde h$ to $M$ gives a Riemannian metric $h:=\iota^*(\tilde 
h\rvert_{\iota(M)})$. Denoting its induced metric by $d^h$ we obtain a metric isometry $\iota\colon (M,d^h)\rightarrow 
(\tilde M, d^{\tilde h})$ and $\iota(M)$ is open and proper --- giving the second and third point. Let $p,q\in M$ with 
$p<q$, i.e., there exists a future directed causal curve $\gamma$ from $p$ to $q$. As $\iota$ is an isometry of $(M,g)$ and 
$(\tilde M, \tilde g)$, the curve $\iota\circ\gamma$ is future directed causal and connects $\iota(p)$ with $\iota(q)$. Thus 
$\iota(p)\tilde <\iota(q)$. The case for $p\ll q$ is completely analogous, giving the fourth point. Finally, let 
$\gamma\colon I\rightarrow M$ be a (locally Lipschitz continuous) curve. Then $\gamma$ is $g$-timelike/causal if and only if 
$\iota\circ\gamma$ is $\tilde g$-timelike/causal by the isometric embedding property. Moreover, by \cite[Prop.\ 2.32]{KS:18} 
we have
\begin{equation*}
 L_\tau(\gamma) = L_g(\gamma) = L_{\tilde g}(\iota\circ\gamma) = L_{\tilde \tau}(\iota\circ\gamma)\,.
\end{equation*}
This gives the fifth point and finishes the proof.
\end{pr}

To illustrate that one can have extensions that are not manifolds we consider the following example, which is a Lorentzian 
version of \cite[Ex.\ 4.2.5]{BBI:01}.
\begin{ex}
Let $\R^2_1$ be two-dimensional Minkowski space and embed it into $\R^3$ as a plane through the origin 
orthogonal to the $z$-direction, i.e., $N:=\{(t,x,0): (t,x)\in\R^2\}$. We now add a half-ray to the origin and give the resulting 
space the structure of a \LLSn. Let $\Gamma:=\{(0,0,z):z\geq 0\}$ and set $\tilde M:= N \cup \Gamma$. On $N$ we use the 
relations from Minkowski space and on $\Gamma$ we define $Z_1:=(0,0,z_1)\ll Z_2:=(0,0, z_2)$ if $z_1 < z_2$, and 
$Z_1\leq Z_2$ if $Z_1\ll Z_2$ or $Z_1=Z_2$. For $p=(t,x,0)\in N$ and $Z\in \Gamma$ we define $p\ll Z$ if $(t,x)\ll 0$ in 
$\R^2_1$ and analogously for the causal relation. 
\begin{figure}[h!]
	\begin{center}
		\includegraphics[width=72mm, height= 48mm]{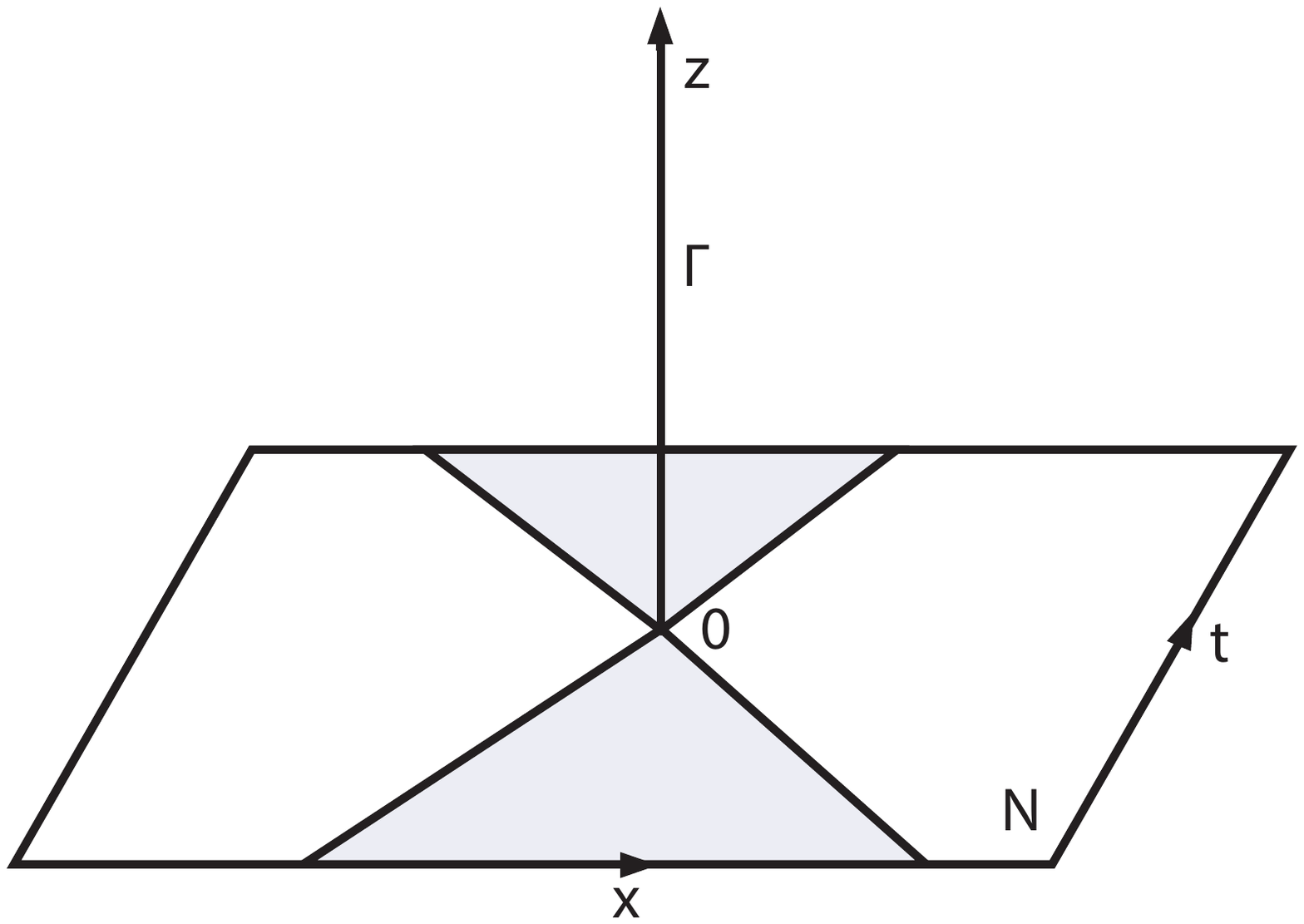}
	\end{center}
	\caption{Non-manifold extension}
\end{figure}
We define the time separation function $\tau$ as the time separation 
function coming from Minkowski space on $N$, for points on $\Gamma$ we set $\tau((0,0,z_1),(0,0,z_2)):=z_2-z_1$ if $z_1\leq 
z_2$ (zero otherwise) and for $p=(t,x,0)\in N$ and $Z=(0,0,z)$ we set $\tau(p,Z):=\sqrt{t^2-x^2} + z$ if $p\leq Z$ (and zero 
otherwise). As $\tau$ is continuous this gives a \LpLSn. In fact, this construction gives a \LLS as it is clearly 
path-connected and locally causally closed. 
Moreover, it is regularly localisable since maximal causal curves always exist 
(they are the, possibly broken, straight lines) and the induced length agrees with the $\tau$-length by construction. 
Furthermore, it is not hard to see that $\tilde M$ is strongly causal. In this space maximal curves branch: every maximal 
curve from $J^-(0)$ to $J^+(0)$ has $0$ as a branching point, as the curve is allowed to continue into $N$ or $\Gamma$. 
This implies via \cite[Cor.\ 4.13]{KS:18} that $\tilde M$ has timelike curvature unbounded below, i.e., a curvature 
singularity in the sense of \cite[Def.\ 4.20]{KS:18}. Finally, $\tilde M$ extends $M\backslash\{(0,0)\}$, thereby providing 
an example of a non-manifold extension. Note that $\tilde M$ does not extend $M$ since $M$ is not embedded into $\tilde M$ as 
an open subset.
\end{ex}

At this point we can introduce the past and future boundary of Lorentzian pre-length spaces with respect to an extension in complete analogy to the 
case of spacetime extensions, see \cite[Def.\ 2.1]{GL:17}.%
\begin{defi}
  Let \tXll be a \LpLS extending the \LpLS \Xll via the embedding $\iota$. The \emph{future/past boundary} 
$\partial^+(X)$\,/\,$\partial^-(X)$ of $X$ is defined as the set of all points $\tilde p\in \partial\iota(X)$ that can be 
reached by a future/past directed $\tilde\ll$-timelike curve $\gamma\colon[0,1]\rightarrow \tilde X$ such that 
$\gamma([0,1))\subseteq \iota(X)$ and $\gamma(1)=\tilde p$.
\end{defi}

The following result establishes that for any extension of a \LLS the future or past boundary is non-empty.
It is a direct analogue of \cite[Lemma 2.17]{Sbi:18}).
\begin{lem}\label{lem:Sbi}
 Let \tXll be an extension of \Xlln, where both are \LLSn s, and denote the corresponding isometry by 
$\iota$. Then there is a $\tilde\ll$-timelike curve $\tilde\gamma\colon[0,1]\rightarrow \tilde X$ such 
that $\tilde\gamma([0,1))\subseteq \iota(X)$ and $\tilde\gamma(1)\in\tilde X\backslash \iota(X)$, i.e., 
$\partial^+(X)\cup\partial^-(X)\neq\emptyset$.
\end{lem}
\begin{pr}
 Since $\iota(X)$ is a proper and open subset of $\tilde X$ and $\tilde X$ is connected, we get that $\partial \iota(X)\neq 
\emptyset$. Let $\tilde p\in \partial \iota(X)$ and let $\tilde \Omega$ be a localising neighbourhood of $\tilde p$ in $\tilde X$. Then, 
$\tilde I^\pm(\tilde p)\cap\tilde\Omega\neq\emptyset$ and let $\tilde q\in \tilde I^-(\tilde p)\cap \tilde \Omega$. We now 
consider two cases. First, if $\tilde q\in\iota(X)$, then since $\tilde q\tilde{\ll}\tilde p$ there is a 
$\tilde{\ll}$-timelike curve $\tilde\gamma\colon[0,1]\rightarrow \tilde X$ such that $\tilde\gamma(0)=\tilde q$, 
$\tilde\gamma(1)=\tilde p$. Set $s_0:=\sup\{s\in[0,1]: \tilde\gamma([0,s])\subseteq\iota(X)\}$, then since $\iota(X)$ is 
open and $\tilde p\in\partial\iota(X)$ we have $\tilde\gamma(s_0)\in\tilde X\backslash\iota(X)$. Reparametrising 
$\tilde\gamma\rvert_{[0,s_0]}$ to $[0,1]$ yields the result. The second case is when $\tilde q\in\tilde 
X\backslash\iota(X)$. Now $\tilde{I}^+(\tilde q)\cap\tilde\Omega$ is a neighbourhood of $\tilde p\in\partial\iota(X)$, thus 
$\iota(X)\cap (\tilde{I}^+(\tilde q)\cap\tilde\Omega)\neq \emptyset$. Let $\tilde r\in \iota(X)\cap (\tilde{I}^+(\tilde 
q)\cap\tilde\Omega)$, then $\tilde q\tilde{\ll}\tilde r$ and the result follows as in the first case by arguing into the 
past.
\end{pr}

\section{Geodesics}\label{sec:geodesics}

In this synthetic approach we have the tools at hand to define causal geodesics as locally length maximising curves. 
Furthermore, we establish that for spacetimes the synthetic notion is compatible with the analytical one.
 
\begin{defi}\label{def:geo}
 Let \Xll be a localising \LpLS and let $\gamma\colon I\rightarrow X$ be a future-directed causal curve. Then $\gamma$ is a 
\emph{geodesic} if for every $t_0\in I$ there exists a localising neighbourhood $\Omega$ of $\gamma(t_0)$ and a
neighbourhood $J=[c,d]$ of $t_0$ in $I$ such that $\gamma\rvert_{J}$ is maximal in $\Omega$ from $\gamma(c)$ to $\gamma(d)$.
\end{defi}

\begin{rem}
 Let $\gamma\colon I\rightarrow X$ be a geodesic and let $t_0\in I$, and $\Omega$ a localising neighbourhood of $\gamma(t_0)$ as 
above. Then
\begin{equation*}
 L_\tau(\gamma\rvert_{[c,d]}) = \omega_{\gamma(t_0)}(\gamma(c),\gamma(d))\,,
\end{equation*}
where $\omega_{\gamma(t_0)}$ is the local time separation function on $\Omega$, cf.\ \cite[Def.\ 3.16]{KS:18}.
\end{rem}

To show that for a smooth and strongly causal spacetime this notion is equivalent to the notion of causal pregeodesics 
we need the following lemma stating a general property of strongly causal Lorentzian length spaces.
\begin{lem}\label{lem-str-cau-tau-om}
 Let \Xll be a strongly causal \LLSn. Then for all $x\in X$ and every localising neighbourhood $\Omega$ of $x$ with local 
time separation function $\omega$ there is a neighbourhood $U$ of $x$, $U\subseteq \Omega$ such that $\omega\rvert_{U\times 
U}$ is completely determined by $\tau$: $\forall p,q\in U:\ \omega(p,q)=\tau(p,q)$. In particular, $\tau$ is 
continuous on a neighbourhood of the diagonal in $X\times X$.
\end{lem}
\begin{pr}
 Let $x\in X$ and let $\Omega$ be a localising neighbourhood of $x$ with local time separation function $\omega$. By strong 
causality and \cite[Lemma 2.38(iii)]{KS:18} there is a neighbourhood $U$ of $x$ with $U\subseteq \Omega$ such that all 
causal curves with endpoints in $U$ are contained in $\Omega$. Let $p,q\in U$ with $p<q$, then by the properties of $\Omega$ 
(see Section \ref{sec:lls}) there is a causal curve $\gamma_{pq}$ that is maximal in $\Omega$ from $p$ to $q$ with 
$L_\tau(\gamma_{pq}) = \omega(p,q)$. As $p,q\in U$, any causal curve connecting these points is contained in $\Omega$. Thus 
$\gamma_{pq}$ is maximal even in $X$, and consequently we have $\tau(p,q) = \Tau(p,q) = L_\tau(\gamma_{pq}) = \omega(p,q)$. 
The neighbourhood of the diagonal can be chosen to be the union of all such $U\times U$ as above.
\end{pr}

With the above lemma we can now establish the promised compatibility.

\begin{thm}\label{thm-pregeo}
Let $(M,g)$ be a smooth, strongly causal spacetime and let $(M,d^h,\ll,\leq,\tau)$ be the induced \LLS (\cite[Ex.\ 
3.24(i)]{KS:18}). Then a causal pregeodesic of $(M,g)$ is a geodesic in the sense of Definition \ref{def:geo} and vice versa.
\end{thm}
\begin{pr}
 First, let $\gamma\colon I\rightarrow M$ be a causal pregeodesic of $(M,g)$, which we can assume without loss of 
generality to be already parametrised as a geodesic. The localising neighbourhoods can be chosen to be (totally) normal 
neighbourhoods. Let $t_0\in I$ and let $U$ be a totally normal neighbourhood of $\gamma(t_0)$. Let $J=[c,d]$ be a
neighbourhood of $t_0$ in $I$ such that 
$\gamma(J)\subseteq U$ and set $x:=\gamma(c), y:=\gamma(d)$. Since $\gamma$ is a geodesic, it 
has to be the radial geodesic from $x$ to $y$ in $U$. As such it is maximal in $U$ and because $L_g=L_\tau$ by \cite[Prop.\ 
2.32]{KS:18} we obtain%
\begin{equation*}
 L_\tau(\gamma\rvert_{[c,d]}) = L_g(\gamma\rvert_{[c,d]}) = 
\sqrt{-g_x(\exp_x^{-1}(y),\exp_x^{-1}(y))}=\omega(x,y)\,.
\end{equation*}

Conversely, let $\gamma\colon I\rightarrow M$ be a geodesic in the sense of Definition \ref{def:geo}. As this is a local 
question, we can cover $\gamma(I)$ by open sets $U$, where $U\subseteq \Omega$ are as 
in the proof of Lemma \ref{lem-str-cau-tau-om}, and show that the 
segment of $\gamma$ in any such $U$ is a pregeodesic with respect to $g$. In fact, let $t_0\in I$ with $\gamma(t_0)\in 
U_0\subseteq \Omega_0$ and let $J\subseteq I$ be an interval around $t_0$ such that $\gamma(J)\subseteq U_0$. Let $s_1, s_2 
\in J$ with $s_1<s_2$, then we get from Lemma \ref{lem-str-cau-tau-om} that
\begin{equation*}
L_\tau(\gamma\rvert_{[s_1,s_2]}) = \omega(\gamma(s_1),\gamma(s_2)) = \tau(\gamma(s_1),\gamma(s_2))\,.
\end{equation*}
Therefore, again since $L_g=L_\tau$, $\gamma$ is maximal on $[s_1,s_2]$ and hence $\gamma$ is a pregeodesic (see e.g.\ \cite[Thm.\ 4.13]{BEE:96}).
\end{pr}
Note that the above proof also shows that the property of being timelike agrees for causal pregeodesics of $(M,g)$ and 
geodesics in the sense of Definition \ref{def:geo} (contrary to the case for arbitrary curves, cf.\ \cite[Ex. 2.22]{KS:18}).

\begin{defi}
  Let \Xll be a localising \LpLS and let $\gamma\colon [a,b)\rightarrow X$ be a future-directed geodesic. Then $\gamma$ 
is \emph{extendible as a geodesic} if there exists a (future-directed) geodesic $\bar{\gamma}\colon[a,b]\rightarrow X$ with 
$\bar{\gamma}\rvert_{[a,b)}=\gamma$. Otherwise, $\gamma$ is called \emph{inextendible as a geodesic}.
\end{defi}

A well-known property of geodesics in smooth semi-Riemannian manifolds is the fact that extendibility
as a geodesic is equivalent to continuous extendibility. Its standard proof relies on the existence of 
convex neighbourhoods. The following result is an analogue in the setting of Lorentzian pre-length spaces,
with localising neighbourhoods working as a substitute. 

\begin{prop} Let \Xll be a strongly causal and localising \LpLS and let $\gamma\colon [a,b)\rightarrow X$ 
be a future-directed geodesic.
Then $\gamma$ is extendible as a geodesic if and only if it is extendible as a continuous curve to $[a,b]$.
\end{prop}
\begin{pr}
Only the `if' part requires a proof, so let us suppose that $\gamma\colon[a,b]\rightarrow X$ is continuous and that 
$\gamma|_{[a,b)}$ is a geodesic. Let $\Omega$ be a localising neighbourhood of $\gamma(b)$ and choose $c\in (a,b)$ such that 
$\gamma([c,b]) \subseteq \Omega$. Then for any $t\in (c,b)$ we have 
\[
L_\tau(\gamma|_{[c,t]}) = \omega(\gamma(c),\gamma(t)),
\]
where $\omega\equiv\omega_{\gamma(b)}$ is the local time separation function on $\Omega$. As $t\nearrow b$,
the right hand side of this equation converges to $\omega(\gamma(c),\gamma(b))$. Concerning the left hand side,
for any $n\in \N$ with $\frac{1}{n}<b-c$ denote by $\gamma_n \colon [c,b]\to X$ a linear reparametrisation of 
$\gamma|_{[c,b-\frac{1}{n}]}$. Then the $\gamma_n$ converge uniformly to $\gamma$ on $[c,b]$. Therefore,
\cite[Prop.\ 3.17]{KS:18} implies that 
\[
L_\tau(\gamma|_{[c,b]}) \ge \limsup_n L_\tau(\gamma_n) = \limsup_n \omega(\gamma(c),\gamma(b-1/n)) = \omega(\gamma(c),\gamma(b)).
\]
As the converse of this inequality holds by the definition of localisability (cf.\ Section \ref{sec:lls}), the claim 
follows.
\end{pr}

\section{Timelike completeness and inextendibility}\label{sec:TC}

As discussed in the introduction, our approach allows us to mimic the proof from the smooth case
that geodesic completeness implies inextendibility, i.e., \cite[Prop.\ 6.16]{BEE:96}. We first
introduce an appropriate notion of timelike geodesic completeness  for \LpLSn s.

\begin{defi}
 Let \Xll be a localising \LpLSn, then $X$ is said to have property $(TC)$ if all inextendible timelike geodesics have infinite 
$\tau$-length.
\end{defi}

This notion is equivalent to timelike geodesic completeness in the case of smooth and strongly causal spacetimes:
\begin{lem}\label{lem-tc}
 Let $(M,d^h,\ll,\leq,\tau)$ be the \LLS induced by a smooth and strongly causal spacetime $(M,g)$. Then $(M,g)$ is  
timelike geodesically complete if and only if $(M,d^h,\ll,\leq,\tau)$ has property $(TC)$.
\end{lem}
\begin{pr}
First, let $(M,g)$ be not timelike geodesically complete, so that there exists an inextendible timelike geodesic 
(without loss of generality inextendible to the future) $\gamma\colon[a,b)\rightarrow M$, with $b<\infty$, thus 
$L_g(\gamma)<\infty$. Since $L_g=L_\tau$ by \cite[Prop.\ 2.32]{KS:18}, Theorem \ref{thm-pregeo} implies that property 
(TC) cannot hold. Conversely, let $(M,g)$ be timelike geodesically complete and let $\gamma\colon[0,b)\rightarrow M$ be an inextendible timelike
geodesic (in the sense of Definition \ref{def:geo}). Then by Theorem \ref{thm-pregeo} $\gamma$ is a timelike pregeodesic of 
$(M,g)$, hence by completeness $L_g(\gamma)=\infty$ (cf.\ \cite[p.\ 154]{ONe:83}). Since $L_g = L_\tau$, property $(TC)$ 
follows.
\end{pr}

Property (TC) does guarantee inextendibility, as the following result shows.

\begin{thm}\label{thm-inext-lls-lls}
 Let \Xll be a strongly causal \LLS that has property $(TC)$. Then \Xll is inextendible as a regular \LLSn.
\end{thm}
\begin{pr}
 Assume, to the contrary, that there exists a regular \LLS \tXll that extends \Xlln. 
 By Lemma \ref{lem:Sbi} there is a (without loss of generality) future-directed $\tilde{\ll}$-timelike curve\\
$\tilde\gamma\colon[0,1]\rightarrow \tilde X$ with $\tilde\gamma([0,1))\subseteq\iota(X)$ and $\tilde\gamma(1)=\tilde p\in 
\tilde X\setminus\iota(X)$. Let $\tilde U$ be a localising neighbourhood of $\tilde p$ (with respect to $\tilde X$) and 
$\tilde \omega$ its local time separation function. Let $t_0\in[0,1)$ be such that $\tilde\gamma([t_0,1])\subseteq \tilde U$. 
Consequently, $q:=\tilde\gamma(t_0)\in \tilde U \cap\iota(X)$ and $q\tilde{\ll} \tilde p$. Thus there is an -- in $\tilde U$ 
-- $\tilde{\tau}$-maximal curve $\tilde\gamma_{q,\tilde p}\colon[0,1]\rightarrow \tilde U$ from $q$ to $\tilde p$, which is 
$\tilde\ll$-timelike by regularity, see \cite[Thm.\ 3.18]{KS:18}. Since 
$\iota(X)$ is open, $q\in\iota(X)$ and $\tilde p\notin\iota(X)$ there is a $t_*\in (0,1)$ such that 
$\tilde\gamma_{q,\tilde p}([0,t_*))\subseteq\iota(X)$ and $\tilde r:=\tilde\gamma_{q,\tilde p}(t_*)\notin\iota(X)$. Then
$\tilde\gamma_{q,\tilde p}\rvert_{[0,t_*)}\colon[0,t_*)\rightarrow \tilde U\cap\iota(X)$ and we set 
$\lambda:=\iota^{-1}\circ \tilde\gamma_{q,\tilde p}\rvert_{[0,t_*)}$. By Definition 
\ref{def-ext},\ref{def-ext-tau}, $\lambda$ is $\ll$-timelike. We claim that $\lambda$ is a timelike $\tau$-geodesic. To 
this end, recall that a maximal causal curve is maximal on any subinterval, see \cite[Prop.\ 2.34,(ii)]{KS:18}. Fix any 
$0\leq s_0 < t_*$, and let $V$ be a neighbourhood of $\lambda(s_0)$ with $\iota(V)\subseteq \tilde U$.
As $X$ is strongly causal, there exists a neighbourhood $V'\subseteq V$ of $\lambda(s_0)$ such that 
any causal curve that starts and ends in $V'$ is contained in $V$.
Now suppose that $s_1\le s_0<s_2$ are so close that $\lambda|_{[s_1,s_2]}$ is contained in $V'$. 
Then in particular any future-directed $\leq$-causal curve connecting 
$\lambda(s_1)$ to $\lambda(s_2)$ remains entirely in $V$. By Definition \ref{def-ext},\ref{def-ext-tau} we 
therefore obtain
\begin{equation*}
\begin{split}
L_\tau(\lambda|_{[s_1,s_2]}) &= L_{\tilde \tau}(\iota\circ\lambda|_{[s_1,s_2]})\\ 
&=\max\{L_{\tilde \tau}(\tilde\alpha):\tilde\alpha \text{ f.d. }  {\tilde\le\text{-causal from }} \iota\circ\lambda(s_1) 
\text{ to } 
\iota\circ\lambda(s_2) \text{ in } \tilde U \}\\
&\ge\max\{L_{\tilde\tau}(\iota\circ\alpha):\alpha \text{ f.d. } {\le\text{-causal from }} \lambda(s_1) \text{ to } 
\lambda(s_2) 
\text{ in } V \}\\
&=\max\{L_{\tau}(\alpha):\alpha \text{ f.d. } {\le\text{-causal from }} \lambda(s_1) \text{ to } \lambda(s_2) \text{ in } V 
\}\\
&=\max\{L_{\tau}(\alpha):\alpha  \text{ f.d. } {\le\text{-causal from }} \lambda(s_1) \text{ to } \lambda(s_2) \text{ in } 
X \}\\
&= \Tau(\lambda(s_1),\lambda(s_2)))\ge L_\tau(\lambda|_{[s_1,s_2]})\,.
\end{split}
\end{equation*}
Thus $L_\tau(\lambda|_{[s_1,s_2]}) = \Tau(\lambda(s_1),\lambda(s_2))) = \tau(\lambda(s_1),\lambda(s_2)))$. By Lemma 
\ref{lem-str-cau-tau-om}, any local time separation function is completely determined by $\tau$ on $V'$, hence the above shows 
that $\lambda$ is a geodesic in $X$. Moreover, the length of $\lambda$ is given by
\begin{equation*}
 L_\tau(\lambda) = L_{\tilde\tau}(\iota\circ\lambda) = \lim_{t\nearrow t_*}L_{\tilde\tau}(\tilde\gamma_{q,\tilde p}\rvert_{[0,t]}) = 
\lim_{t\nearrow t_*}\tilde\omega (q,\tilde\gamma_{q,\tilde p}(t)) = \tilde\omega(q,r)<\infty\,,
\end{equation*}
as the local time separation function $\tilde\omega$ of $\tilde U$ (with respect to $\tilde X$) is continuous and finite. 
Finally, $\lambda$ is inextendible as a geodesic in $X$ since it is not even extendible as a continuous curve ($\lim_{t\nearrow 
t_*}\iota\circ\lambda(t) = \lim_{t\nearrow t_*}\gamma_{q,\tilde p}(t)= \tilde r\notin\iota(X)$) --- thus contradicting property 
$(TC)$.
\end{pr}

We can now relate the low regularity inextendibility to a blow-up of curvature. More precisely, we have 
the following result.

\begin{thm}\label{thm-tc-ext-cur-sing}
 Let \Xll be a strongly causal \LLS that has property $(TC)$. Suppose that $X$ can be extended
to a strongly causal locally timelike geodesically connected \LLSn. 
Then this extension has a causal curvature 
singularity (\cite[Def.\ 4.20]{KS:18}). Specifically, the extension cannot have bounded upper causal curvature.
\end{thm}
\begin{pr}
Let \Xll be a \LLS that is strongly causal and has property $(TC)$. Assume that there exists a \LLS \tXll extending \Xll satisfying the assumptions of the theorem and having causal curvature 
bounded above. Then \cite[Rem.\ 4.16, Thm.\ 4.17 and Thm.\ 4.18]{KS:18} yield that \tXll is regular. 
This contradicts the inextendibility result Theorem \ref{thm-inext-lls-lls} and yields that $X$ has a curvature singularity 
in the sense of \cite[Def.\ 4.20]{KS:18}.
\end{pr}

We now specialise to the case where the object to be extended is a smooth spacetime. Firstly, recall that~\emph{causally 
plain\/} spacetimes are precisely those that do not exhibit the bubbling phenomenon. Roughly speaking, a metric is bubbling 
if it contains a point where the boundary of the future null cone has non-empty interior. (For a precise definition, 
see~\cite[Definition~1.16]{CG:12}; cf.~also the discussion preceding Lemma~5.6 in~\cite{KS:18}.) Spacetimes $(M, g)$ with 
$g$ a Lipschitz metric are causally plain~\cite[Corollary~1.17]{CG:12}. 

The following result is now a direct corollary of Theorem \ref{thm-inext-lls-lls}.

\begin{cor}\label{cor-geo-compl-inext-lls}
 Let $(M,g)$ be a smooth, strongly causal and timelike geodesically complete spacetime and let $(M,d^h,\ll,\leq,\tau)$ be its 
induced \LLSn. Then $(M,d^h,\ll,\leq,\tau)$ is inextendible as a regular Lorentzian length space, and hence also inextendible in the class of continuous, 
strongly causal and causally plain spacetimes that are regular.
\end{cor}
\begin{pr}
 By Lemma \ref{lem-tc} $(M,d^h,\ll,\leq,\tau)$ has property $(TC)$ and strong causality is the same notion for  
spacetimes and the corresponding \LLSn s by \cite[Lemma 2.21(i),(ii) and Lemma 2.38(iii)]{KS:18}. 
Thus, Theorem \ref{thm-inext-lls-lls} applies, showing that 
$(M,d^h,\ll,\leq,\tau)$ is inextendible as a regular \LLSn. Furthermore, by \cite[Thm.\ 5.12]{KS:18} every continuous
strongly causal and causally plain spacetime $(\tilde M, \tilde g)$ gives rise to a \LLSn.
\end{pr}

Also in this case of spacetimes we obtain the result that timelike geodesic completeness forces the extension to have a 
curvature singularity, even though curvature cannot be defined in the usual sense via the Riemann 
tensor.

\begin{cor}\label{cor-inext-st-cur-sing}
  Let $(M,g)$ be a smooth, strongly causal and timelike geodesically complete spacetime and let $(M,d^h,\ll,\leq,\tau)$ be 
its induced \LLSn. If $(M,d^h,\ll,\leq,\tau)$ is extendible as a strongly causal locally
timelike geodesically connected \LLS then the extension has a causal curvature singularity 
(it cannot have causal curvature bounded above).
\end{cor}
\begin{pr}
This follows directly from Theorem \ref{thm-tc-ext-cur-sing}, similarly to the proof of Corollary \ref{cor-geo-compl-inext-lls}.
\end{pr}

\begin{rem} 
In \cite{AB:08}, Alexander and Bishop introduced sectional curvature bounds for general semi-Riemannian manifolds. Moreover, 
they characterized these curvature bounds via triangle comparison with small triangles in model spaces (i.e., the spaces 
$M_K$ from \eqref{eq:model_spaces} in the Lorentzian setting), see \cite[Thm.\ 1.1]{AB:08}. As was shown in \cite[Ex.\ 
4.9]{KS:18}, our definitions in 
Section \ref{sec:lls} are compatible with these curvature bounds in this sense and in particular a curvature 
singularity in our sense implies that there cannot be a corresponding sectional curvature 
bound in the sense of \cite{AB:08}. 
Corollary \ref{cor-inext-st-cur-sing} therefore implies that if the extension is assumed to be a smooth and strongly causal 
spacetime itself, then its sectional curvature as defined in \cite{AB:08} must be unbounded above.
\end{rem}
To conclude this section we note that it is an interesting open question whether one can characterize completeness of 
timelike geodesics in \LLSn s via condition $(TC)$, analogous to the smooth case, cf.\ \cite[p.\ 154]{ONe:83}. 
\section{Relation to other results on low regularity inextendibility}\label{sec:relation}
In this final section we relate our work to further current results on the low regularity inextendibility of spacetimes.

In \cite{GL:18} it was recently established that in a (locally) Lipschitz continuous spacetime maximal causal curves have a 
causal character. This immediately gives that the induced \LLS $(M,d^h,\ll,\leq,\tau)$ of a strongly causal Lipschitz 
spacetime $(M,g)$ is regular: By \cite[Cor.\ 1.17]{CG:12} and \cite[Thm.\ 5.12]{KS:18} $(M,d^h,\ll,\leq,\tau)$ is a \LLS and 
by \cite[Thm.\ 1.1]{GL:18} it is regular (a fact that was already observed by Graf and Ling in \cite{GL:18}). From this they 
deduce that a timelike geodesically complete smooth spacetime is inextendible in the class of Lipschitz spacetimes. Thus, 
their result is slightly stronger than ours when restricted to spacetimes (compare Corollary \ref{cor-geo-compl-inext-lls}) 
as they do not need strong causality of the original spacetime. However, even when restricting to the case where the object 
to be extended is a spacetime, our result is more general  in the following sense:
\begin{itemize}
 \item It allows the original spacetime to be of low regularity (continuous and causally plain) as well.
 \item There might be continuous strongly causal, causally plain spacetimes inducing a regular \LLS where the metric is not locally 
Lipschitz continuous.
 \item It applies even to non-manifold extensions, and
 \item it relates inextendibility with curvature blow-up (Theorem \ref{thm-tc-ext-cur-sing}).
\end{itemize}

In \cite{GLS:18} the authors show that a smooth, timelike geodesically complete and globally hyperbolic spacetime is 
$\Con^0$-inextendible, i.e., there is no spacetime with continuous metric extending the given spacetime. Again, as above, 
their result is slightly stronger when restricting to spacetimes, since of course not all spacetimes with continuous metrics 
give rise to a \LLSn, as they need not be causally plain and strongly causal (see e.g.\ \cite[Ex.\ 1.11]{CG:12}). However, 
our approach does not need the original spacetime to be globally hyperbolic and (as above) allows it to be of low regularity 
as well. Moreover, as noted above our result also rules out non-manifold extensions (as long as they are regular \LLSn s). A 
closer inspection of the proof of Theorem \ref{thm-inext-lls-lls} reveals that one does not need that the entire extension is 
regular. In fact, all that is needed is that a maximal causal curve $\gamma$ that is contained in the original space except for
its endpoint (which is on the boundary) is timelike whenever its starting point and endpoint are timelike related in the 
extension. This is weaker than being regular, as it essentially only concerns points in the original space and its boundary. 
Thus the main result of \cite{GLS:18} can be understood in this way: If the smooth spacetime is timelike geodesically complete 
and globally hyperbolic, then maximal causal curves as above have a causal character. This then yields the inextendibility 
result.

It should also be noted that in our framework one can define \emph{future/past one-connectedness} (\cite[Def.\ 
2.13]{Sbi:18}) and \emph{future/past divergence} (\cite[Def.\ 2.4(2)]{GL:17}) as for spacetimes. Since being extendible forces 
the future or past boundary to be non-empty by Lemma \ref{lem:Sbi} a further line of study could be to see if, as for 
spacetimes, future (past) one-connectedness together with future (past) divergence yield empty future (past) boundary (cf.\ 
\cite[Thm.\ 2.5]{GL:17}).
\bigskip

To summarize, we have developed a framework where we can show inextendibility of spaces that resemble timelike 
geodesically complete spacetimes, in a similar spirit as the classical result (\cite[Prop.\ 6.16]{BEE:96}). Our approach provides a 
new and unified perspective on the recent results \cite{GLS:18, GL:18}, see the discussion above. 
Moreover, for the first time we can relate low regularity inextendibility with a (synthetic) curvature blow up --- a fact 
that fits well with physical expectations. Finally, it shows that timelike geodesic completeness is a very robust property, 
which carries over even to spaces that are not spacetimes or even manifolds.

\bigskip

\noindent
{\bf Acknowledgements.} This work was supported by research grants P26859 and P28770 of the Austrian Science Fund FWF. The 
work of J.G.\ was partially supported by STFC Consolidated Grant ST/L000490/1.


\addcontentsline{toc}{section}{References}

\end{document}